\documentclass[letterpaper]{article} % DO NOT CHANGE THIS
\usepackage{aaai23}  % DO NOT CHANGE THIS
\usepackage{times}  % DO NOT CHANGE THIS
\usepackage{helvet}  % DO NOT CHANGE THIS
\usepackage{courier}  % DO NOT CHANGE THIS
\usepackage[hyphens]{url}  % DO NOT CHANGE THIS
\usepackage{graphicx} % DO NOT CHANGE THIS
\usepackage{natbib}  % DO NOT CHANGE THIS AND DO NOT ADD ANY OPTIONS TO IT
\usepackage{caption} % DO NOT CHANGE THIS AND DO NOT ADD ANY OPTIONS TO IT
\usepackage{algorithm}
\usepackage{algorithmic}

\usepackage{newfloat}
\usepackage{listings}
\DeclareCaptionStyle{ruled}{labelfont=normalfont,labelsep=colon,strut=off} % DO NOT CHANGE THIS
\floatstyle{ruled}
\newfloat{listing}{tb}{lst}{}
\floatname{listing}{Listing}
\usepackage{bm}
\usepackage{multirow}
\usepackage{subfigure}
\usepackage{amsmath}
\usepackage{amssymb}
\nocopyright %-- Your paper will not be published if you use this command

\usepackage[colorlinks,linkcolor=red]{hyperref}
\title{DMIS: Dynamic Mesh-based Importance Sampling for Training Physics-Informed Neural Networks}
\author{
    Zijiang Yang, \textsuperscript{\rm 1,2}
    Zhongwei Qiu, \textsuperscript{\rm 1,2,3}
    Dongmei Fu \textsuperscript{\rm 1,2,4}
}
\affiliations{
    \textsuperscript{\rm 1}School of Automation and Electrical Engineering, University of Science and Technology Beijing\\
    \textsuperscript{\rm 2}Beijing Engineering Research Center of Industrial Spectrum Imaging\\
    \textsuperscript{\rm 3}The University of Sydney\\
    \textsuperscript{\rm 4}Shunde Innovation School, University of Science and Technology Beijing\\

    zijiangyang@xs.ustb.edu.cn, qiuzhongwei@xs.ustb.edu.cn,  fdm\_ustb@ustb.edu.cn
}

\usepackage{bibentry}

\begin{document}

\maketitle

\begin{abstract}

Modeling dynamics in the form of partial differential equations (PDEs) is an effectual way to understand real-world physics processes.
For complex physics systems, analytical solutions are not available and numerical solutions are widely-used.
However, traditional numerical algorithms are computationally expensive and challenging in handling multiphysics systems.
Recently, using neural networks to solve PDEs has made significant progress, called physics-informed neural networks (PINNs).
PINNs encode physical laws into neural networks and learn the continuous solutions of PDEs.
For the training of PINNs, existing methods suffer from the problems of inefficiency and unstable convergence, since the PDE residuals require calculating automatic differentiation.
In this paper, we propose \textbf{D}ynamic \textbf{M}esh-based \textbf{I}mportance \textbf{S}ampling (DMIS) to tackle these problems.
DMIS is a novel sampling scheme based on importance sampling, which constructs a dynamic triangular mesh to estimate sample weights efficiently.
DMIS has broad applicability and can be easily integrated into existing methods.
The evaluation of DMIS on three widely-used benchmarks shows that DMIS improves the convergence speed and accuracy in the meantime.
Especially in solving the highly nonlinear Schrödinger Equation, compared with state-of-the-art methods, DMIS shows up to 46\% smaller root mean square error and five times faster convergence speed.
Code is available at https://github.com/MatrixBrain/DMIS.

\end{abstract}

\section{Introduction}

Modeling the dynamics of real-world physics systems has important guiding significance for production activities, such as fluid mechanics and heat transfer.
These physics systems are usually described by partial differential equations (PDEs). 
Due to the highly nonlinear of these PDEs, in most cases, analytical solutions are not available.

\begin{figure}[!ht]
    \centering
    \includegraphics[width=0.9\columnwidth]{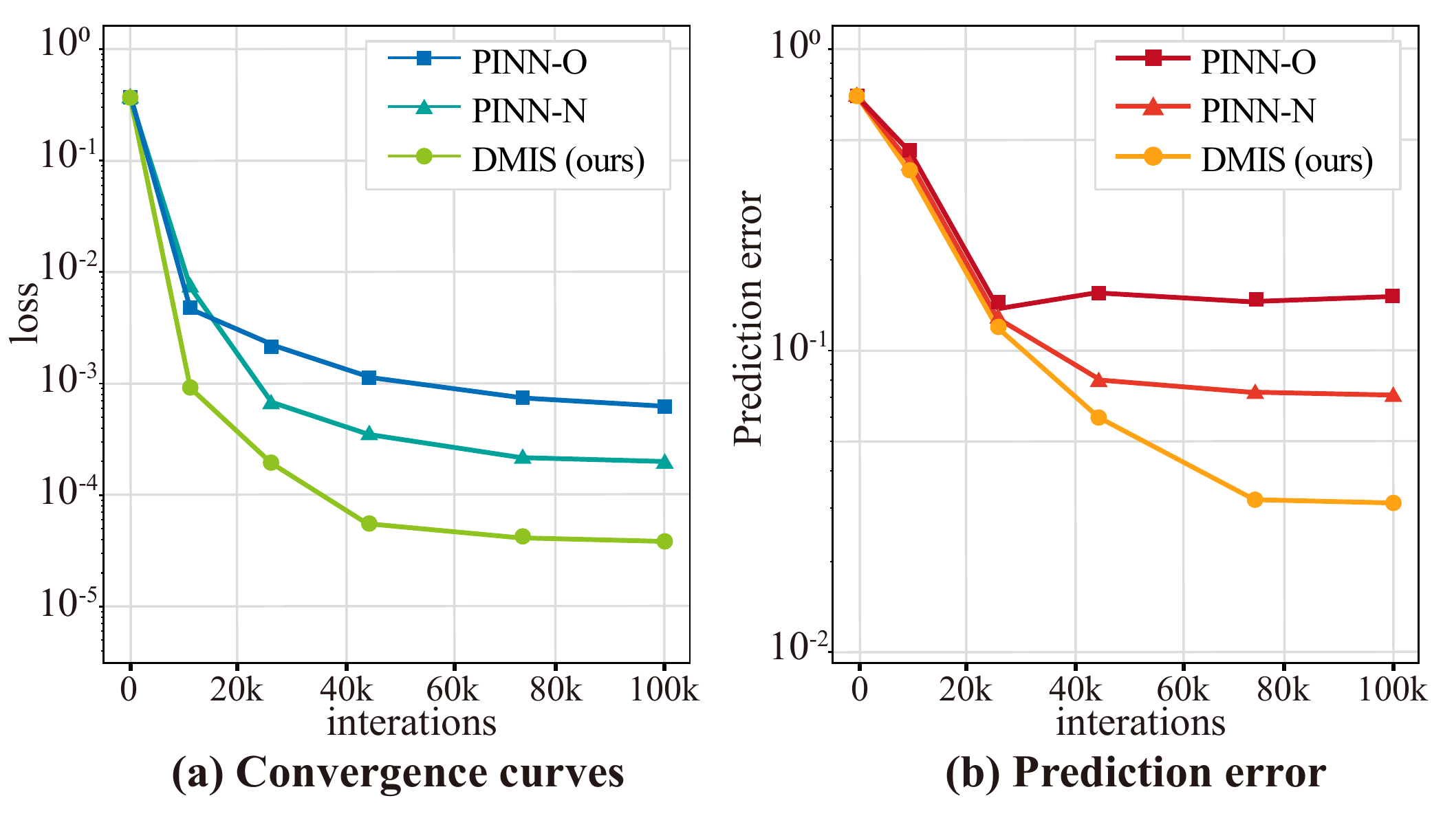}
    \vspace{-0.4cm}
    \caption{
    Comparison with stat-of-the-art methods (PINN-O, PINN-N). DMIS achieves faster convergence speed and better prediction accuracy.
    }
    \label{fig:1stFig}
\end{figure}

In the past decades, numerical algorithms for solving PDEs have been greatly developed.
Although classical numerical algorithms significantly promoted the development of related fields, these algorithms are computationally expensive and face severe challenges in multiphysics and multiscale systems.
With the improvement of data acquisition capability, it is an important issue to use data to modify simulation results.
In general, data assimilation is the primary method to combine classical numerical algorithms and observation data, but data assimilation introduces additional uncertainties and seriously affects convergence. 

Recently, with the rapid development of machine learning, the ability to extract features and mine information from observation data has been dramatically improved. \citep{krizhevsky2012imagenet, he2016deep}.
As an effective supplement to classical numerical methods, Physics-informed Neural Networks (PINNs) formulate the problem of solving PDEs into a parameter optimization problem \citep{raissi2019physics}.
PINNs encode PDEs into the loss function of neural networks.
In addition, boundary conditions and initial conditions are also integrated into the loss function as soft constraints.
PINNs are mesh-free and learn continuous solutions of PDEs.
Practices show that PINNs can be applied to solve different types of PDEs\citep{raissi2019physics, pang2019fpinns, zhang2020learning}.
Since partial derivative terms of PDE residuals require automatic differentiation to calculate, the training of PINNs is computationally expensive and unstable.
Especially for solving PDEs with high order partial derivative terms, these defects are more prominent.
Therefore, in the studies of PINNs, it is a critical issue to design schemes to stabilize training and improve training efficiency.

Three mainstream improvement directions are deeply studied for the training of PINNs, including weights of loss terms, parallel computation, and sampling.
The loss function of PINNs contains multiple loss terms, and methods of weight allocation are mainly to balance these terms.
These methods usually need to be customized according to physics systems, which limits the applications.
The efficient parallel computation of PINNs can be realized by decomposing the whole domain into several subdomains.
Nevertheless, domain decomposition involves additional subdomain boundary conditions and affects the convergence efficiency and model accuracy.
The last mainstream improvement direction, sampling, improves training efficiency and model accuracy by changing the sampling probability of data. Due to broad applicability, it attracts increasing attention.

The training of PINNs involves two sampling processes. 
The first sampling process is to build a training dataset from the whole domain, and the second is to sample mini-batch data in each iteration for parameter optimization.
In order to distinguish, in this paper, we refer to the first sampling process as generating collocation points.
Currently, the progress in sampling mainly focuses on the generation of collocation points.
These collocation points generation approaches improve the training efficiency by constructing a better training dataset without incurring additional computational costs.
The mini-batch sampling significantly impacts convergence speed and model accuracy, which has been verified in other fields \citep{shrivastava2016training, kang2019decoupling, mildenhall2020nerf, lei2021less}.
However, the mini-batch sampling for PINNs has not been fully studied.
In uniform sampling, the training spends most of computation on data points that are not helpful for optimization.

In our previous research, we found that training efficiency and model accuracy can be improved by introducing well-designed weighting algorithm in other fields\citep{qiu2019learning, qiu2020dgcn, wang2021learning, qiu2022dynamic}.
In this paper, we introduce additional weighting algorithm and propose a novel sampling scheme, called Dynamic Mesh-based Importance Sampling (DMIS), to speed up the training of PINNs.
To guarantee the sampling method theoretically, we introduce the concept of importance sampling into DMIS.
However, importance sampling requires calculating the sampling probability of each point, which leads to high computational costs.
To reduce the computational cost, we propose a novel sampling weight estimation method, called Dynamic Mesh-based weight estimation (DMWE), which constructs a dynamic triangular mesh to estimate the weight of each data point efficiently.
The mesh constructed by DMWE is updated dynamically according to the loss distribution in the whole domain during training.
DMIS efficiently integrates importance sampling into the training of PINNs, and has a low computational cost and broad applicability.
As shown in Figure \ref{fig:1stFig}, DMIS achieves faster convergence speed and better accuracy compared with the state-of-the-art methods.

Our contributions can be summarized as follows:
\begin{itemize}
    \item We propose an efficient importance sampling scheme for training PINNs, which improves the training convergence speed and the model accuracy without significantly increasing computational cost.
    \item We propose a method to calculate sample weights in PINNs efficiently. In addition to be used in our importance sampling scheme, this method is also suitable for other approaches that need to calculate sample weights.
    \item The extensive experiments on three widely-used benchmarks demonstrate the superior performance of our importance sampling scheme.
\end{itemize}

\section{Related Work}

\subsubsection{Physics-Informed Neural Networks}
In recent years, benefited from the significant progress of deep learning, the method of solving PDEs based on neural networks has been dramatically developed.
According to the method of embedding physics laws, all of these methods based on neural networks for solving PDEs can be grouped into three types, including observational bias \citep{umetani2018learning}, inductive bias \citep{cai2020phase, dong2021method} and learning bias \citep{sirignano2018dgm, raissi2019physics}.
Among these methods, physics-informed neural networks attract increasing attention, owing to excellent augment ability \citep{niaki2021physics, raissi2020hidden, wandel2022spline}.

The training of PINNs is computationally expensive and unstable.
To alleviate these defects of PINNs, allocation of loss weights \citep{wang2021understanding, wang2022and, krishnapriyan2021characterizing}, parallel computation methods \citep{meng2020ppinn, karniadakis2020extended} and sampling schemes \citep{lu2021deepxde, nabian2021efficient, daw2022rethinking, wu2022comprehensive} have been widely discussed and improve the convergence efficiency and prediction accuracy of PINNs. 
Gradient control methods are also studied \cite{kim2021dpm}.

\begin{figure*}[!ht]
    \centering
    \includegraphics[width=0.90\textwidth]{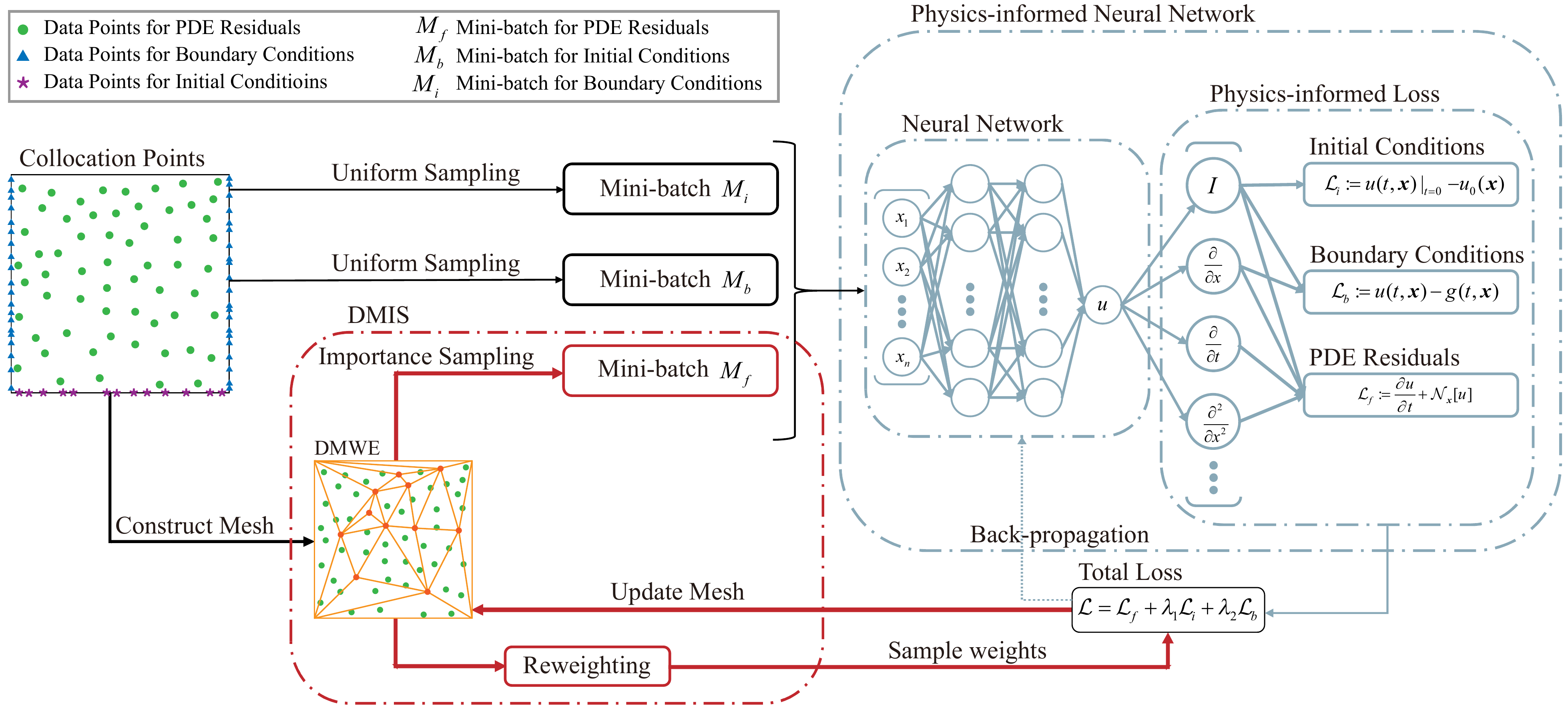}
    \caption{
    The Training of PINNs with DMIS.
    In each iteration, sample weights are estimated by DMWE, and a mini-batch for PDE residuals is sampled by DMIS.
    The derivatives terms are computed by automatic differentiation.
    These mini-batches are applied to compute the initial conditions loss $\mathcal{L}_i$, boundary conditions loss $\mathcal{L}_b$, and PDE residuals loss $\mathcal{L}_f$.
    DMIS reweights sample weights of $M_f$ and conveys these weights to the loss function $\mathcal{L}$.
    The mesh of DMIS is updated according to $\mathcal{L}$.
    }
    \label{fig:Diagram}
\end{figure*}

\subsubsection{Importance Sampling}
Importance sampling is an effective method to accelerate Monte Carlo Integration, which has been applied to train neural networks \citep{schaul2016prioritized, chen2018fastgcn, meng2022count}.
It has been proved theoretically that the SGD has the fastest convergence speed if the mini-batch is sampled according to a unique distribution that the sampling probability of each point is proportional to the 2-norm of loss gradient with respect to parameters  \citep{needell2014stochastic, zhao2015stochastic}.
The computational cost of the theoretical optimal method is prohibitive. Follow-up works mainly focus on finding approximation methods \citep{loshchilov2015online, canevet2016importance, johnson2018training, katharopoulos2018not}.
To our best knowledge, only \citet{nabian2021efficient} preliminary explores the combination of importance sampling and the training of PINNs.

\section{Optimization Problem of PINNs}

PINNs learn a approximate solution $\hat{u}(t, \bm{x}; \bm{\theta})$ to fit the latent solution $u(t, \bm{x})$ of the following PDE:
\begin{equation}
    \begin{split}
        &\frac{\partial{u}}{\partial{t}} + \mathcal{N}_{\bm{x}} [u] = 0, \bm{x} \in \Omega, t\in[0, T], \\
        &u(t, \bm{x})|_{t=0} = u_0(\bm{x}), \bm{x} \in \Omega, \\
        &u(t, \bm{x})=g(t, \bm{x}), x\in \partial{\Omega}, t\in[0, T],
    \end{split}
    \label{Eq:ProblemDefine}
\end{equation}
where $\bm{\theta}$ is the parameter of PINNs, $\mathcal{N}_{\bm{x}}$ denotes a differential operator consisted of spatial derivatives, $u_{0}(\bm{x})$ is the initial condition, $g(\bm{x})$ is the boundary condition, $\bm{x}$ is a D-dimensional position vector, and $\Omega$ is a subset of $\mathbb{R}^D$ with boundary $\partial{\Omega}$.
For the convenience of subsequent discussion, the input vector composed of time $t$ and space vector $\bm{x}$ is denoted as $\mathrm{\bm{x}}$.
The optimization goal of PINNs is to minimize the residual of PDEs with the constraints of satisfying boundary conditions and initial conditions:
\begin{equation}
    \begin{split}
        &\bm{\theta^*} = \arg \min_{\bm{\theta}} r_f(\bm{\theta}),\\
        &s.t.~ r_i(\bm{\theta})=0, r_b(\bm{\theta}) = 0,
    \end{split}
    \label{Eq:OpProblemRaw}
\end{equation}
where $r_f(\bm{\theta})$, $r_i({\bm{\theta}})$ and $r_b({\bm{\theta}})$ are the residuals of PDEs, initial conditions, and boundary conditions, respectively. 

Equation (\ref{Eq:OpProblemRaw}) is difficult to solve.
PINNs regard constraints as penalty terms and formulate the constrained optimization problem into an unconstrained optimization problem:
\begin{equation}
    \bm{\theta^*} = \arg \min_{\bm{\theta}} r_f(\bm{\theta}) + \lambda_1 r_i(\bm{\theta}) + \lambda_2 r_b(\bm{\theta}),
    \label{Eq:OpProblemPINNs}
\end{equation}
where $\lambda_1$ and $\lambda_2$ are weights.
The common practice of PINNs is to fit $\bm{\theta}$ by Monte Carlo approximation.
PINNs generate collocation points from spatial-temporal domains uniformly, and mini-batch stochastic gradient descent method (SGD) is employed to optimize parameters.
Focusing on data points helpful for parameter optimization is a more efficient sampling strategy.
Monte-Carlo approximation has provided the mathematical tools, called importance sampling, to design such a sampling method.

\section{Method}

Sampling approaches significantly impact training efficiency.
Recent works mainly focus on the generation of collocation points, while the mini-batch sampling is neglected.
Motivated by the theoretical completeness of importance sampling and the non-negligible impact of mini-batch sampling, we design a novel sampling approach for mini-batch sampling based on importance sampling to improve the convergence speed and model accuracy of PINNs.

\subsubsection{Importance Sampling for PINNs}
PINNs generate collocation points from domain and domain boundary. Equation (\ref{Eq:OpProblemPINNs}) can be approximated as a loss function of data points:
\begin{equation}
    \bm{\theta}^* 
    \approx \arg\min_{\bm{\theta}}\mathcal{L}_f + \lambda_1\mathcal{L}_i + \lambda_2 \mathcal{L}_b,
\label{Eq:OpProblemReal}
\end{equation}
where $\mathcal{L}_f$, $\mathcal{L}_i$, and $\mathcal{L}_b$ are the losses of PDE residuals, initial conditions, and boundary conditions, respectively. 

Datasets for PDE residuals, initial conditions and boundary conditions are denoted by $N_f$, $N_i$ and $N_b$, respectively.
In general, mini-batches are uniformly sampled from $N_f$, $N_i$, and $N_b$, respectively.
According to the Monte Carlo approximation, we can introduce a more efficient sampling method.
Because boundary conditions and initial conditions are penalty terms, we only introduce importance sampling into the sampling of $N_f$.
The loss of PDE residuals $\mathcal{L}_f$ combined with importance sampling is shown as:
\begin{equation}
    \mathcal{L}_f = \frac{1}{\vert N_f \vert}\sum_{i=1}^{|N_f|}\alpha_i\ell_f(\mathrm{x}_i;\bm{\theta}), ~\alpha_i = \frac{p_{i}}{q_{i}},
\end{equation}
where $\vert N_f \vert$ is the size of $N_f$, $\ell_{f}$ is the PDE residual of each data point, $\alpha_i$, $p_{i}$ and $q_{i}$ are sample weight, sampling probability and alternative sampling probability of data point $\mathrm{x}_i$, respectively.
Considering mini-batches are obtained by uniform sampling in general, $p_{i}$ equals to $1/\vert N_f \vert$ for $i\in\{1, 2, \cdots, \vert N_f\vert \}$ and the calculation of $\alpha_i$ is shown as:

\begin{equation}
    \alpha_i = \frac{1}{\vert N_f\vert q_i}, ~ i\in\{1, 2, \cdots, |N_f|\},
    \label{Eq:SampleWeightInit}
\end{equation}

\subsubsection{Simplified Calculation \& Reweighting} 

The critical issue of importance sampling is to determine $q_i$ for $i \in \{ 1, 2, \cdots, \vert N_f \vert\}$ and we hope to find the best alternative sampling distribution to make the convergence rate fastest. Suppose convergence rate is defined as:
\begin{equation}
\begin{split}
    C^{(t)} = - \mathbb{E}_f[\Vert \bm{\theta}^{(t+1)} - \bm{\theta}^{*} \Vert_2^2 - \Vert \bm{\theta}^{(t)} - \bm{\theta}^{*} \Vert_2^2],
\end{split}
    \label{Eq:ConvergenceRateDefine}
\end{equation}
where $C^{(t)}$ is the convergence rate at step $t$, $\bm{\theta}^{(t)}$ and $\bm{\theta}^{(t+1)}$ are parameters at step $t$ and $t+1$, respectively.

With the definition of convergence rate shown in Equation (\ref{Eq:ConvergenceRateDefine}), \citet{needell2014stochastic, zhao2015stochastic} have demonstrated that the best sampling probability of collocation points is determined by $ q^* \propto \Vert \nabla_{\bm{\theta}} \ell_f(\mathrm{x}, \bm{\theta})\Vert_2 $.

However, the computational cost of this theoretical optimal sampling method is unacceptable, and it is necessary to find alternative methods.
Inspired by \citet{katharopoulos2017biased}, an approximation calculation of the theoretical optimal formula is employed in DMIS.
\begin{equation}
    q_i^{(t)} = \frac{\ell_f(\mathrm{x}_i, \bm{\theta}^{(t)})}{\sum_{j=1}^{N_f}\ell_f(\mathrm{x}_j, \bm{\theta}^{(t)})}, ~ i\in\{1, 2, \cdots, |N_f|\},
    \label{Eq:SimplifiedFormula}
\end{equation}
where $q_i^{(t)}$ is the sampling probability of $\mathrm{x}_i$ at step $t$.
\citet{katharopoulos2017biased} also demonstrate that Equation (\ref{Eq:SimplifiedFormula}) does not change the rank of sampling probability. $\forall i, j \in \{1, 2, \cdots, \vert N_f \vert\}$, if $q_i^* < q_j^*$, we can obtain that $q_i^{(t)} < q_j^{(t)}$.
Therefore, Equation (\ref{Eq:SimplifiedFormula}) is a reasonable approximation. However, we find that sample weights calculated by Equation (\ref{Eq:SimplifiedFormula}) lead to unstable training in the initial stage. This problem is caused by data points with high loss, which lead to sharp local gradients. To fix this problem, we introduce a super parameter $\beta$ to adjust $\alpha$. With $\beta > 1$, greater penalties are applied to data points with high loss and the result is denoted as $\alpha^\prime$:
\begin{equation}
    \alpha_i^\prime = (\frac{1}{\vert N_f\vert q_i})^\beta, ~ \beta\in[1, +\infty), ~ i\in\{1, 2, \cdots, |N_f|\}.
    \label{Eq:Reweighting}
\end{equation}

\begin{table*}[!ht]
    \centering
    \begin{tabular*}{\hsize}{@{}@{\extracolsep{\fill}}c c c c c c c c c c@{}}
    \hline
    
    \hline
    \multicolumn{1}{c}{\multirow{2}{*}{PDE}} & \multicolumn{2}{c}{Network Config} & Optimizer Config & \multicolumn{3}{c}{DMIS Config} & \multicolumn{3}{c}{Dataset Config}\\ \cline{2-10} & Depth & Width & Learning rate & $\vert S \vert$ & $\gamma$ & $\beta$ & $\vert N_{f} \vert $ & $\vert N_{i} \vert $ & $\vert N_{b} \vert $ \\ \hline
        Schrödinger & 4 & 64 & 0.001 & 1000 & 0.4 & 2 & 60000 & 200 & 200 \\ 
        Burgers & 3 & 32 & 0.005 & 1000 & 0.4 & 1.5 & 100000 & 2000 & 2000\\
        KdV & 4 & 64 & 0.001 & 1000 & 0.4 & 2 & 60000 & 2000 & 2000\\
        \hline
        
        \hline
    \end{tabular*}
    \caption{The hyper-parameters used for each benchmark. $\vert S \vert$ is the set size of mesh points, $\gamma$ is the mesh update threshold, and $\beta$ is the hyper-parameter of reweighting. $\vert N_f \vert$, $\vert N_i \vert$, and $\vert N_b \vert$ are the dataset size of PDE residuals, initial conditions and boundary conditions, respectively.}
    \label{Tb:hyperparameters}
\end{table*}
\begin{table*}[!ht]
    \centering
    \begin{tabular*}{\hsize}{@{}@{\extracolsep{\fill}}cccccccccc@{}} 
    \hline
    
    \hline
    \multicolumn{1}{c}{\multirow{2}{*}{Method}} &
    \multicolumn{3}{c}{Schrödinger Equation} &
    \multicolumn{3}{c}{Burgers' Equation} &
    \multicolumn{3}{c}{KdV Equation} \\
    \cline{2-10} & ME & MAE & RMSE & ME & MAE & RMSE & ME & MAE & RMSE \\ \hline
    PINN-O & 1.360 & 0.186 & 0.4092 & 0.451 & 0.0738 & 0.1100 & 2.140 & 0.363 & 0.520 \\ 
    PINN-N & 0.948 & 0.149 & 0.2906 & 0.358 & 0.0579 & 0.0859 & 1.860 & 0.292 & 0.441 \\
    xPINN & 0.546 & 0.045 & 0.0089 & 0.261 & 0.0099 & 0.0010 & 2.462 & 0.272 & 0.230 \\
    cPINN & 0.591 & 0.069 & 0.0169 & 0.324 & \textbf{0.0084} & \textbf{0.0007} & 2.925 & 0.258 & 0.248 \\
    \hline
    \textbf{PINN-DMIS(ours)} & 0.647 & 0.127 & 0.2196 & \textbf{0.225} & 0.0294 & 0.0495 & \textbf{1.170} & 0.391 & 0.492 \\
    \textbf{xPINN-DMIS(ours)} & 0.867 & 0.036 & 0.0129 & 0.420 & 0.0115 & 0.0017 & 2.380 & 0.233 & \textbf{0.196} \\
    \textbf{cPINN-DMIS(ours)} & \textbf{0.358} & \textbf{0.025} & \textbf{0.0033} & 0.397 & 0.0111 & 0.0016 & 2.680 & \textbf{0.230} & 0.200 \\
    \hline
    
    \hline
    \end{tabular*}
    \caption{Comparison with PINN-O \citep{raissi2019physics}, PINN-N \citep{nabian2021efficient}, xPINN \citep{karniadakis2020extended} and cPINN \citep{jagtap2020conservative} on benchmarks of the Schrödinger Equation, the Viscous Burgers' Equation and the KdV Equation.}
    \label{Tb:AccEvaluation}
\end{table*}

\begin{algorithm}[tb]
\caption{Sampler with DMIS}
\label{alg:Sampler}
\textbf{Input}: batch size of PDE residuals $\vert M_f \vert$. \\
\textbf{Parameter}: set size of mesh points $\vert S \vert$, reweighting parameter $\beta$, mesh update threshold $\gamma$, dataset of PDE residuals $N_f$, iteration step $t$. \\
\textbf{Output}: mini-batch $M_f$, vector of sample weights $\bm{\alpha^{\prime}}$. \\
\textbf{Initialization}
\begin{algorithmic}[1]
\STATE $t_0 \leftarrow 0$
\STATE $q^{(t_0)}_{i} \leftarrow 1/\vert N_f \vert, i\in\{1, 2, \cdots, \vert N_f \vert\}$
\STATE $g_{i}^{(t_0)} \leftarrow 1/|N_f|, i\in\{1, 2, \cdots, \vert N_f \vert\}$
\STATE $S \leftarrow \vert S \vert $ points sampled with $ g_{i}^{(t_0)}$ from $N_f$
\STATE Build triangular mesh by Delaunay Triangulation
\end{algorithmic}

\textbf{Mini-batch Sampling}
\begin{algorithmic}[1]

\STATE Compute $\{ \ell_f(\mathrm{x}_i, \bm{\theta}^{(t-1)})\}_{\mathrm{x}_i\in S}$
\STATE Estimate score of other points in $N_f$ by interpolation
\STATE Compute $q_{i}^{(t)}$ according to Equation (\ref{Eq:SimplifiedFormula})
\STATE $M_f \leftarrow \vert M_f \vert$ points sampled with $q_{i}^{(t)}$ from $N_f$
\STATE Compute $\alpha_{i}^{(t)}$ according to Equation (\ref{Eq:SampleWeightInit})
\STATE Compute $\alpha_{i}^{\prime(t)}$ according to Equation (\ref{Eq:Reweighting})
\STATE Compute $\mathrm{Sim}(\bm{v}^{(t_0)}, \bm{v}^{(t)})$ according to Equation (\ref{Eq:CosBasedEvaluation})

\IF{$\mathrm{Sim}(\bm{v}^{(t_0)}, \bm{v}^{(t)}) < \gamma$}
\STATE $g_{i}^{(t)} \propto \Vert q^{(t-1)}_{i} - q^{(t_0)}_{i} \Vert, ~ i \in \{1, 2, \cdots, \vert N_f \vert\}$
\STATE $q_{i}^{(t_0)} \leftarrow q_{i}^{(t)}, ~ i \in \{1, 2, \cdots, \vert N_f \vert\}$
\STATE $S \leftarrow \vert S \vert $ points sampled with $ g_{i}$ from $N_f$
\STATE $t_0 \leftarrow t$
\STATE Update mesh by Delaunay Triangulation
\ENDIF

\STATE \textbf{return} $M_f$, $\bm{\alpha^{\prime}}$
\end{algorithmic}
\end{algorithm}

\subsubsection{DMWE}
Since Equation (\ref{Eq:SimplifiedFormula}) reduces the computational cost of each data point, the sampling probability still needs to be calculated point by point, and it is an enormous burden in solving complex PDEs.
To further reduce the computational cost, we propose dynamic mesh-based weight estimation (DMWE) to calculate sample weight by interpolation.

In DMWE, interpolation based on Delaunay Triangulation is employed. 
Specifically, We dynamically generate a subset $S$ from $N_f$ to construct a triangular mesh.
DMWE only calculated the sample weights of points in $S$ exactly, and the weights of other points are obtained by interpolation.
$S$ is generated according to Equation (\ref{Eq:BasicSeeds}).
\begin{equation}
    g_i^{(t)} \propto \Vert q_i^{(t)} - q_i^{(t-1)}\Vert, ~ i \in \{1, 2, \cdots, \vert N_f \vert\},
    \label{Eq:BasicSeeds}
\end{equation}
where $g_i^{(t)}$ is the selection probability of point $\mathrm{x}_i$ at step $t$.
Equation (\ref{Eq:BasicSeeds}) reduces the number of mesh points in the inactive region to reduce the computational cost. Meanwhile, Equation (\ref{Eq:BasicSeeds}) also ensures high-precision interpolation in the active region.

The interpolation based on Delaunay is time-consuming and it is also unnecessary to update the triangular mesh in each iteration step.
Therefore, we introduce the cosine similarity-based evaluation method to decide whether to re-select $S$ and rebuild the triangular mesh.
\begin{equation}
    \mathrm{Sim}(\bm{v}^{(t_0)}, \bm{v}^{(t)}) = \frac{\bm{v}^{(t_0)} \cdot \bm{v}^{(t)}}{||\bm{v}^{(t_0)}|| \cdot \Vert \bm{v}^{(t)} \Vert},
    \label{Eq:CosBasedEvaluation}
\end{equation}
where $v^{(t_0)}$ and $v^{(t)}$ are vectors composed of the sample weight of data points in $S$ at step $t_0$ and step $t$, respectively.
If the cosine similarity is smaller than the threshold $\gamma$, $S$ will be re-selected from $N_f$, and the mesh will be updated.

\subsubsection{DMIS}
The pseudo-code of a sampler with DMIS is shown as Algorithm \ref{alg:Sampler}.
In each iteration, the PDE residual of data points in $S$ are computed first, and then the sample weights of other points are estimated by interpolation.
The training process combined with DMIS is shown as Figure \ref{fig:Diagram}.

\section{Experiments}

\subsection{Benchmark}

We consider to solve Schrödinger equation, Viscous Burgers' equation and Korteweg-de Vries equation. PINNs have excellent interpolation accuracy, but the extrapolation accuracy still needs to be further improved \citep{kim2021dpm}. 
For this reason, the data division method used in our work is different from that in original PINNs \citep{raissi2019physics}.
Specifically, for each solving problem, we use the similar method as \citet{kim2021dpm} to divide the entire time domain $[0, T]$ into three segments: $[0, T/2]$, $[T/2, 3T/4]$ and $[3T/4, T]$. 
These three segments are used for training, validating, and testing.

\subsubsection{Schrödinger Equation}
Schrödinger equation is a fundamental equation in quantum mechanics.
Schrödinger equation with an initial condition $u(0, x) = 2 sech(x)$ and periodic boundary conditions is considered:
\begin{equation}
i\frac{\partial{u}}{\partial{t}} + 0.5\frac{\partial^2{u}}{\partial{x^2}} + \vert u \vert^2 u = 0, x\in[-5, 5], t\in[0, \pi/2].
\end{equation}

\subsubsection{Viscous Burgers' Equation}
The Burgers' equation simulates shock wave propagation and reflection.
We consider the Burgers' equation with an initial condition $u(0, x) = -\sin(\pi x)$ and a boundary condition $u(t, x) = 0, x\in\{-1, 1\}$, which is shown as follows:
\begin{equation}
\frac{\partial{u}}{\partial{t}} + u\frac{\partial{u}}{\partial{x}}=\frac{0.04}{\pi}\frac{\partial^2{u}}{\partial{x^2}}, x\in[-1, 1], t\in[0, 1].
\end{equation}

\subsubsection{Korteweg-de Vries Equation}
Korteweg-de Vries (KdV) equation describes the waves on shallow water surfaces.
The KdV equation, which involves third-order derivatives, is suitable to evaluate the efficiency of DMIS in solving partial differential equations with higher-order derivatives.
We consider the KdV equation with an initial condition $u(0, x) = cos(\pi x)$ and periodic boundary conditions:
\begin{equation}
\frac{\partial{u}}{\partial{t}} + u\frac{\partial{u}}{\partial{x}} + 0.0025 \frac{\partial^3{u}}{\partial{x^3}} = 0, x\in[-1, 1], t\in[0, 1].
\end{equation}

\subsection{Experimental Setting}

\subsubsection{Baseline}
We evaluate the performance of DMIS based on original PINN\citep{raissi2019physics}, xPINN\citep{karniadakis2020extended} and cPINN\citep{jagtap2020conservative}.
The original PINN is denoted by PINN-O.
PINN with importance sampling scheme of \citet{nabian2021efficient} is denoted by PINN-N.
PINN, xPINN and cPINN with our approach, DMIS, are denoted by PINN-DMIS, xPINN-DMIS and cPINN-DMIS, respectively.

\subsubsection{Hyper-parameter} 
Table \ref{Tb:hyperparameters} summarizes hyper-parameters of networks, optimizers, datasets, and DMIS for each benchmark. 
Collocation points are uniformly generated from the domain and domain boundary.
The Adam optimizer is employed on all benchmarks.
For PINN-N, we uniformly sample 10,000 seeds within the domain. For cPINN and xPINN, the spatial is equally decomposed into three subdomains.

\begin{figure}[!ht]
    \centering
    \includegraphics[width=0.9\columnwidth]{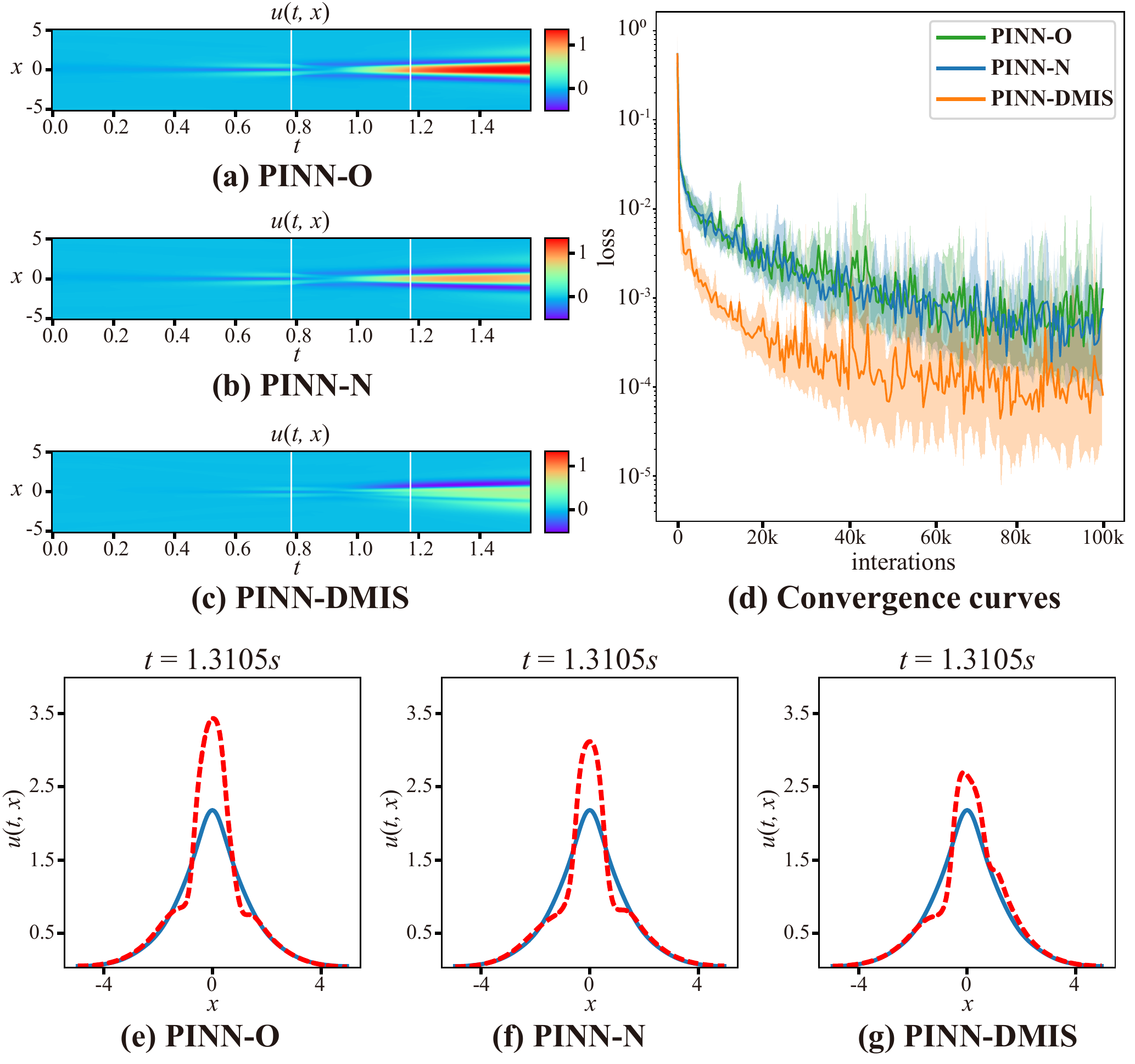}
    \caption{The Schrödinger Equation. (a, b, c) The prediction errors. (d) Convergence curves. (e, f, g) show the prediction (red) and ground-truth (blue).}
    \label{fig:Schrodinger}
\end{figure}

\subsubsection{Evaluation Metrics}
We choose the maximum error (ME), mean absolute error (MAE), and root mean square error (RMSE) as evaluation metrics of prediction accuracy.
The exact numerical solutions are obtained by Py-PDE \citep{zwicker2020py}.
We introduce the calculation time and the number of iteration steps required for convergence, denoted by $TC$ and $NC$, respectively, to evaluate the impact of sampling methods on convergence behavior.
The subscript indicates the convergence level. For example, $NC_{5}$ is the minimum iteration steps required when the training loss is stable below 1e-5 for 1000 iterations.
Considering that the convergence behavior is challenging to measure, $NC$ and $TC$ are only used to evaluate the convergence behavior in this paper roughly.
It is more intuitive to observe convergence behavior through convergence curves.

\begin{table}[!t]
    \centering
    \begin{tabular*}{\hsize}{@{}@{\extracolsep{\fill}}c c c c c@{}}
    \hline
    
    \hline
    ~ & $TC_2/s$ & $TC_3/s$ & $NC_2$ & $NC_3$ \\ \hline
    PINN-O & 151 & 2005 & 4740 & 61984 \\
    PINN-N & 466 & 5968 & 5681 & 72959 \\ \hline
    \textbf{PINN-DMIS(ours)} & \textbf{34} & \textbf{399} & \textbf{847} & \textbf{10219} \\
    \hline
    
    \hline
    \end{tabular*}
    \caption{Evaluation results of convergence behavior for solving the Schrödinger Equation}
    \label{Tb:SchrodingerConvergenceBehavior}
\end{table}

\subsection{Experimental Result}

\subsubsection{Schrödinger Equation}
Figure \ref{fig:Schrodinger}(a, b, c) report the prediction error by PINN-O, PINN-N, and PINN-DMIS for solving Schrödinger Equation.
Compared with other methods, PINN-DMIS has better performance.
The snapshots in Figure \ref{fig:Schrodinger} show that PINN-DMIS has the lowest maximum errors.
Table \ref{Tb:AccEvaluation} summarizes the results of prediction accuracy for solving Schrödinger Equation.
DMIS can significantly improve the accuracy.
Figure \ref{fig:Schrodinger}(d) reports that PINN-DMIS converges fastest and has the lowest training error.
Table \ref{Tb:SchrodingerConvergenceBehavior} compares the convergence behavior.
PINN-DMIS shows up to five times faster convergence speed than PINN-O.

\begin{figure}[!ht]
    \centering
    \includegraphics[width=0.9\columnwidth]{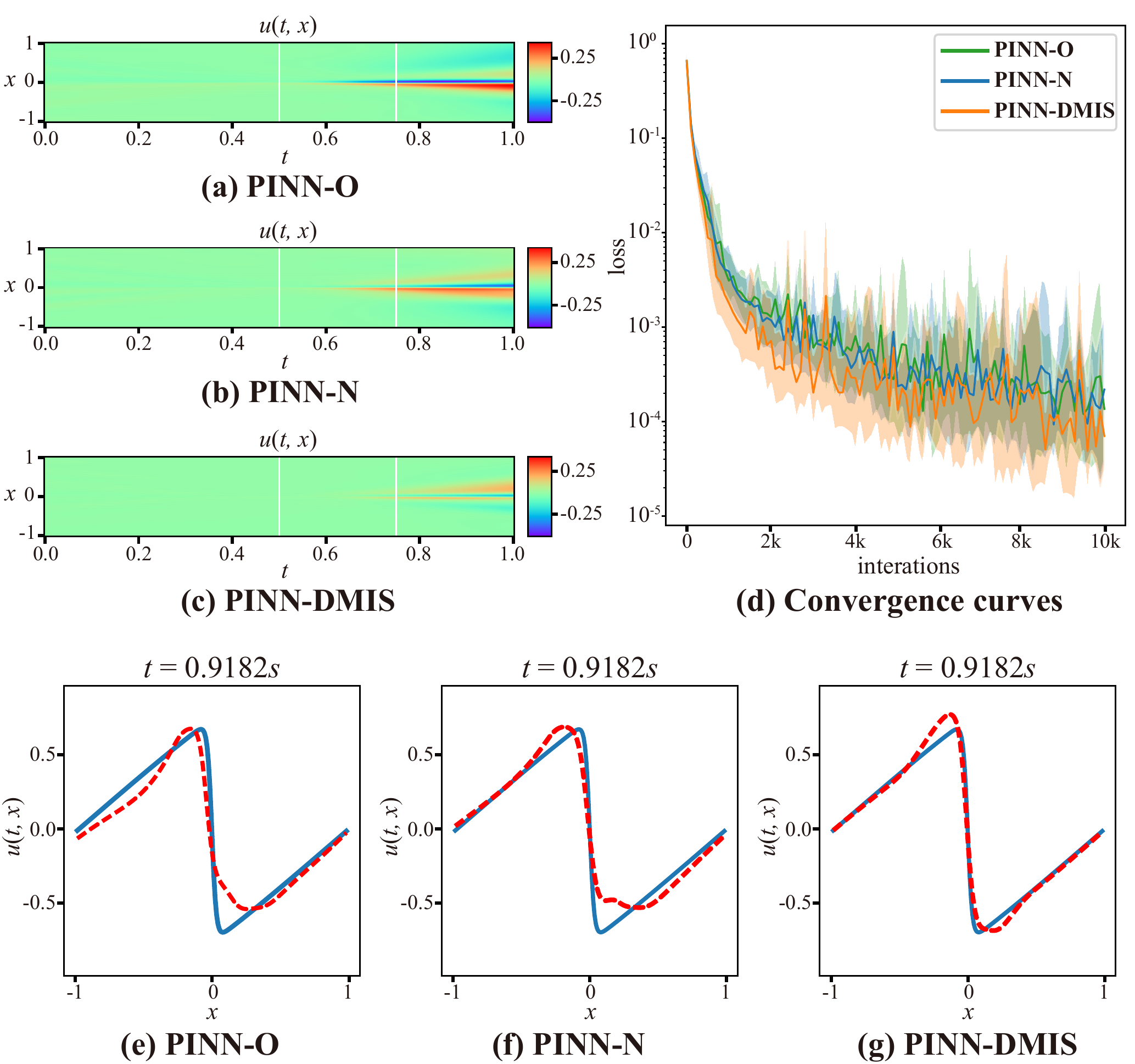}
    \caption{The Burgers' Equation. (a, b, c) The prediction errors. (d) Convergence curves. (e, f, g) show the prediction (red) and ground-truth (blue).}
    \label{fig:Burgers}
\end{figure}

\subsubsection{Viscous Burgers’ Equation}
Figure \ref{fig:Burgers}(a, b, c) report the prediction errors for solving the Burgers’ Equation.
PINN-DMIS has lower prediction error in the extrapolation segment.
Figure \ref{fig:Burgers}(e, f, g) report the snapshots of prediction at t=0.9182s.
PINN-O and PINN-N fail to predict $u(t, x)$, and PINN-DMIS still has perfect prediction performance.
Table \ref{Tb:AccEvaluation} summarizes the evaluation results for solving the Viscous Burgers' Equation.
PINNs-DMIS has the lowest maximum error.
The convergence curves for solving the Viscous Burgers' Equation are reported in Figure \ref{fig:Burgers}(d).
PINN-DMIS has a slender lead convergence speed.
Table \ref{Tb:BurgersConvergenceBehavior} reports the convergence behavior of PINN-O, PINN-N, and PINN-DMIS for solving the Viscous Burgers' Equation.
These three methods have similar convergence behavior in the initial stage, but in the following stage, PINN-DMIS converges fastest.

\begin{table}[!t]
    \centering
    \begin{tabular*}{\hsize}{@{}@{\extracolsep{\fill}}c c c c c@{}}
    \hline
    
    \hline
    ~ & $TC_2/s$ & $TC_3/s$ & $NC_2$ & $TC_3$ \\ \hline
    PINN-O & \textbf{15} & 97 & 1218 & 7308\\
    PINN-N & 83 & 417 & 1218 & 6032\\ \hline
    \textbf{PINN-DMIS(ours)} & 20 & \textbf{89} & \textbf{1160} & \textbf{5336}\\
    \hline
    
    \hline
    \end{tabular*}
    \caption{Evaluation results of convergence behavior for solving the Viscous Burgers' Equation}
    \label{Tb:BurgersConvergenceBehavior}
\end{table}

\begin{table*}[!ht]
    \centering
    \begin{tabular*}{\hsize}{@{}@{\extracolsep{\fill}}cccccccccc@{}} 
    \hline
    
    \hline
    \multicolumn{1}{c}{\multirow{2}{*}{Method}} &
    \multicolumn{3}{c}{Schrödinger Equation} &
    \multicolumn{3}{c}{Burgers' Equation} &
    \multicolumn{3}{c}{KdV Equation} \\
    \cline{2-10} & ME & MAE & RMSE & ME & MAE & RMSE & ME & MAE & RMSE \\ \hline
    PINN-O & 1.360 & 0.186 & 0.409 & 0.451 & 0.074 & 0.110 & 2.140 & \textbf{0.363} & 0.520 \\ \hline
    \textbf{PINN-BasicIS(ours)} & 0.965 & 0.140 & 0.279 & 0.273 & \textbf{0.025} & \textbf{0.044} & 1.460 & 0.376 & 0.526 \\ 
    \textbf{PINN-DMIS(ours)} & \textbf{0.647} & \textbf{0.127} & \textbf{0.220} & \textbf{0.225} & 0.029 & 0.049 & \textbf{1.170} & 0.391 & \textbf{0.492} \\ 
    \hline
    
    \hline
    \end{tabular*}
    \caption{Ablation Study of DMWE for solving the Schrödinger Equation, the Viscous Burgers' Equation and the KdV Equation.}
    \label{Tb:AblationDMWE}
\end{table*}
\begin{table*}[!ht]
    \centering
    \begin{tabular*}{\hsize}{@{}@{\extracolsep{\fill}}ccccccccccc@{}} 
    \hline
    
    \hline
    \multirow{2}{*}{$\vert S \vert$} & 
    \multirow{2}{*}{$\gamma$} &
    % \multicolumn{1}{c}{\multirow{2}{*}{$\vert S \vert$}} &
    % \multicolumn{1}{|c|}{\multirow{2}{*}{$\gamma$}} &
    \multicolumn{3}{c}{Schrödinger Equation} & 
    \multicolumn{3}{c}{Burgers' Equation} & 
    \multicolumn{3}{c}{KdV Equation}\\
    \cline{3-11}
    ~ & ~ & ME & MAE & RMSE & ME & MAE & RMSE & ME & MAE & RMSE \\ \hline
    1000 & 0.2 & 0.762 & 0.120 & 0.233 & 0.593 & 0.039 & 0.085 & 1.460 & 0.376 & 0.526 \\
    1000 & 0.4 & \textbf{0.647} & 0.127 & 0.220 & \textbf{0.225} & \textbf{0.029} & \textbf{0.049} & \textbf{1.170} & 0.391 & 0.492 \\
    1000 & 0.6 & 0.838 & \textbf{0.111} & \textbf{0.208} & 0.919 & 0.057 & 0.147 & 1.760 & \textbf{0.351} & 0.518 \\
    10000 & 0.4 & 0.913 & 0.132 & 0.262 & 0.392 & 0.033 & 0.075 & 1.280 & 0.366 & \textbf{0.475}\\
    20000 & 0.4 & 1.100 & 0.163 & 0.332 & 0.620 & 0.046 & 0.105 & 1.730 & 0.452 & 0.600\\ 
    \hline
    
    \hline
    \end{tabular*}
    \caption{Ablation studies of the set size $\vert S \vert$ and the mesh update threshold $\gamma$ for solving the Schrödinger Equation, the Viscous Burgers' Equation and the KdV Equation.}
    \label{Tb:AblationSuperparameters}
\end{table*}

\subsubsection{Korteweg-de Vries Equation}
\begin{figure}[!ht]
    \centering
    \includegraphics[width=0.9\columnwidth]{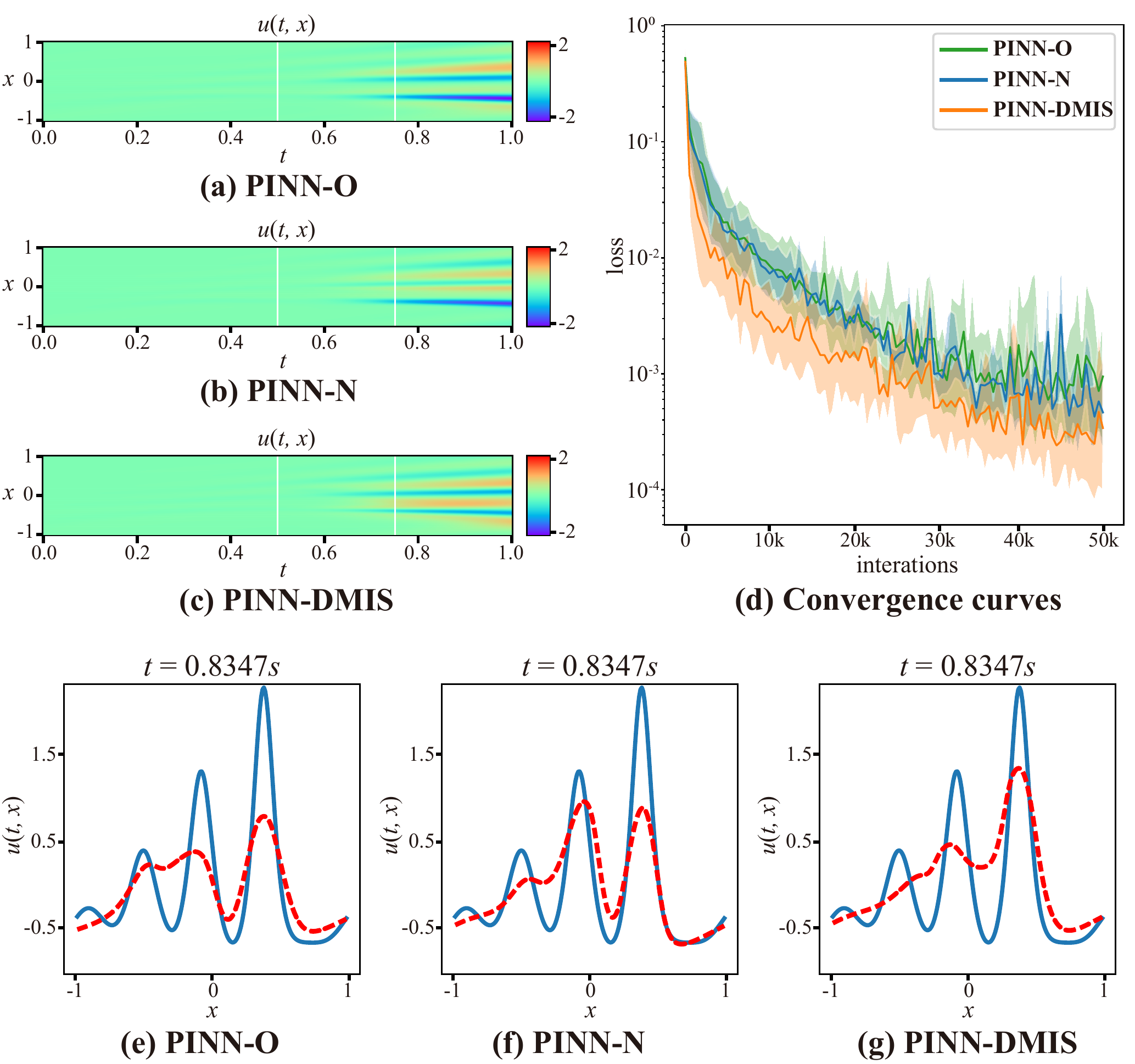}
    \caption{The KdV Equation. (a, b, c) The prediction errors. (d) Convergence curves. (e, f, g) show the prediction (red) and ground-truth (blue).}
    \label{fig:KDV}
\end{figure}

Figure \ref{fig:KDV}(e, f, g) report the snapshots of prediction at t=0.8347s. all methods fail to predict $u(t, x)$.
By contrast, PINN-DMIS has lower maximum error.
Table \ref{Tb:AccEvaluation} reports prediction accuracy and DMIS can effectively improve the accuracy of various PINNs.
Figure \ref{fig:KDV}(d) reports the convergence curves for solving the KdV Equation.
PINN-DMIS converges faster than other methods.
Table \ref{Tb:KdVConvergenceBehavior} reports the evaluation results of convergence behavior for solving the KdV Equation.
The results indicate that PINN with DMIS has noticeable convergence acceleration.

\subsection{Ablation Study}
% To further demonstrate the effectiveness of DMIS, we study the impact of DMWE and hyper-parameters.

\subsubsection{DMWE}
\begin{figure}[!ht]
    \centering
    \includegraphics[width=0.95\columnwidth]{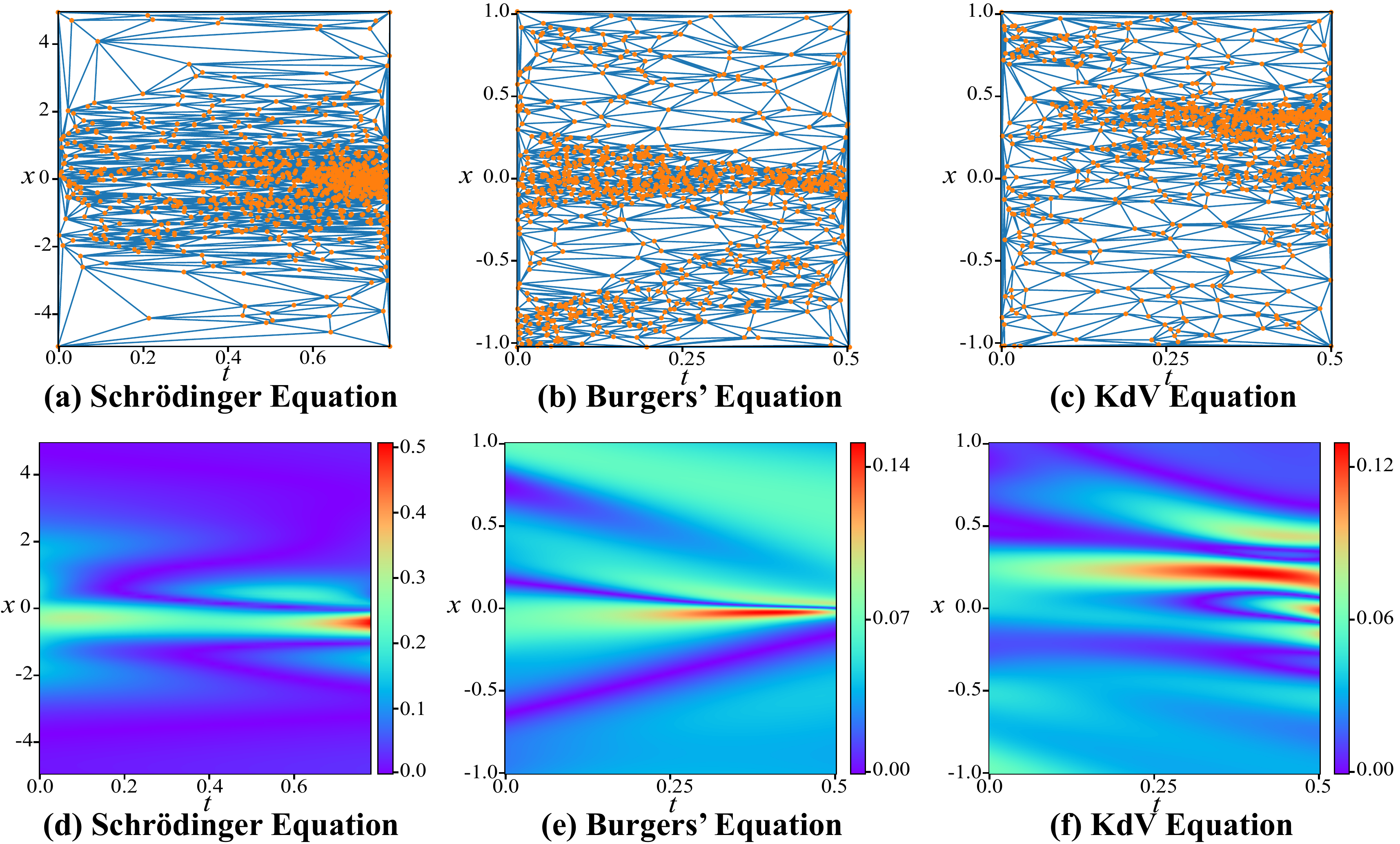}
    \caption{Dynamic mesh. (a, b, c) Meshes constructed by DMWE. (d, e, f) The corresponding changes of $\mathcal{L}$ in adjacent updates.}
    \label{fig:MeshUpdate}
\end{figure}

Figure \ref{fig:MeshUpdate} reports meshes constructed by DMWE and the corresponding changes of $\mathcal{L}$ in adjacent updates.
During training, mesh points gradually gather in the region with a drastic change of $\mathcal{L}$, and samples in these regions significantly impact the parameter optimization in the current stage.
DMWE aims to arrange more mesh points in the active region to achieve a high-precision estimation of sample weight and arrange fewer mesh points in the inactive region to reduce the computational cost.
The dynamic meshes shown in Figure \ref{fig:MeshUpdate} meet our expectations.

Our approach without DMWE is denoted by PINN-BasicIS.
Table \ref{Tb:AblationDMWE} reports evaluation results of PINN-O, PINN-BasicIS, and PINN-DMIS.
Compared with PINN-BasicIS, PINN-DMIS has better performance in most cases.
Especially for solving Schrödinger Equation, PINN-DMIS has the best performance on all evaluation metrics.

\subsubsection{Hyper-parameter}
Table \ref{Tb:AblationSuperparameters} reports the ablation studies of mesh update threshold $\gamma$ and set size $\vert S \vert$.
When $\vert S \vert$ is too tiny, the estimation deviation of sample weights increases dramatically. 
On the other hand, when $\vert S \vert$ is too large, the weights of samples are unstable, which also leads to a decline of model performance.
Similarly, when $\gamma$ is too tiny, DMWE fails to update mesh.
When $\gamma$ is too large, the mesh is frequently updated, leading to a severe decrease of training efficiency.
Fortunately, $\gamma=0.4$ is always appropriate and we recommend setting $\gamma$ to 0.4.

\begin{table}[!t]
    \centering
    \begin{tabular*}{\hsize}{@{}@{\extracolsep{\fill}}c c c c c@{}}
    \hline
    
    \hline
    ~ & $TC_2/s$ & $TC_3/s$ & $NC_2$ & $NC_3$ \\ \hline
    PINN-O & 476 & 1812 & 12631 & 48401 \\
    PINN-N & 794 & 4001 & 9855 & 49912\\ \hline
    \textbf{PINN-DMIS(ours)} & \textbf{288} & \textbf{1046} & \textbf{6954} & \textbf{25029}\\
    \hline
    
    \hline
    \end{tabular*}
    \caption{Evaluation results of convergence behavior for solving the KdV Equation}
    \label{Tb:KdVConvergenceBehavior}
\end{table}

\section{Conclusion}

In this paper, we propose a novel importance sampling scheme, Dynamic Mesh-based Importance Sampling (DMIS).
DMIS constructs a dynamic triangular mesh to estimate sample weights and effectively integrates importance sampling into the training of PINNs.
Experiments on three widely-used benchmarks show that DMIS significantly improves the convergence speed and accuracy of PINNs without significantly increasing computational cost.
We also evaluate the performance of DMIS on three different PINNs models and results verify effectiveness and excellent generalization ability of DMIS.

\clearpage
\section{Acknowledgments}

This work is supported by the National Key Research and Development Program of China (No.2021YFB3702402) and the Scientific and Technological Innovation of Shunde Graduate School of University of Science and
Technology Beijing (No. BK20AE004).

\bibliography{paper}

\clearpage

\appendix

\section{Appendix A: Supporting Theory}
In this section, we supplement some deductions.

\subsection{Exact Importance Sampling for Neural Networks}

Parameters of neural networks are optimized by training. 
The goal of training is as follows.
\begin{equation}
\begin{split}
    \bm{\theta}^{*} 
    &= \arg\min_{\bm{\theta}} \mathcal{L},\\
    &= \arg\min_{\bm{\theta}} \frac{1}{\vert N\vert}\sum_{\mathrm{x}\in N}\ell(\mathrm{x}; \bm{\theta}),
\end{split}
\end{equation}
where $N$ is a dataset, $\mathrm{x}$ is a data point in $N$, $\bm{\theta}$ is parameters of neural networks, $\bm{\theta}^{*}$ is the optimal parameters, $\ell$ is a loss function, and $\mathcal{L}$ is the total loss.

Parameters of neural networks are updated according to the gradient of $\mathcal{L}$:
\begin{equation}
\begin{split}
\bm{\theta}^{(t+1)} 
&= \bm{\theta}^{(t)} - \eta \nabla_{\bm{\theta}^{(t)}} \mathcal{L}, \\
&= \bm{\theta}^{(t)} - \eta \nabla_{\bm{\theta}^{(t)}} \frac{1}{\vert N\vert}\sum_{\mathrm{x}\in N}\alpha_\mathrm{x}\ell(\mathrm{x}; \bm{\theta}),
\end{split}
\label{Eq: AppendixParametersUpdate}
\end{equation}
where $\alpha_\mathrm{x}$ is a sample weight,  $\bm{\theta}^{(t)}$ and $\bm{\theta}^{(t+1)}$ are parameters at iterations $t$ and $t+1$, respectively, and $\eta$ is learning rate.
For uniform sampling, $\alpha_\mathrm{x}=1$. 
Assuming that the sampling probability of $\mathrm{x}$ is $q_{\mathrm{x}}$, with the condition of unbiased estimation, the sample weight $\alpha_\mathrm{x}=1/(\vert N \vert q_{\mathrm{x}})$. 
We aim to find the best sampling probability $q_{\mathrm{x}}$ to make the training convergence rate the fastest.

In our paper, the convergence rate of training is defined as follows \citep{needell2014stochastic, zhao2015stochastic}.
\begin{equation}
\begin{split}
    C^{(t)} = - \mathbb{E}[\Vert \bm{\theta}^{(t+1)} - \bm{\theta}^{*} \Vert_2^2 - \Vert \bm{\theta}^{(t)} - \bm{\theta}^{*} \Vert_2^2],
\end{split}
\label{Eq:AppendixConvergenceRateDefine}
\end{equation}
where $C^{(t)}$ is the convergence rate at iteration $t$, and $\mathbb{E}$ denotes expectation.

For the convenience of subsequent discussion, We set $G^{(t)} = \eta \nabla_{\bm{\theta}^{(t)}} \mathcal{L}$, $D^{(t)} = \Vert \bm{\theta}^{(t)} - \bm{\theta}^{*} \Vert_2^2$, and $D^{(t)} = \Vert \bm{\theta}^{(t)} - \bm{\theta}^{*} \Vert_2^2$.

$D^{(t)}$ and $D^{(t+1)}$ can be simplified as follows. 
\begin{equation}
\small
\begin{split}
    D^{(t)}
    &= (\bm{\theta}^{(t)} - \bm{\theta}^{*})^T(\bm{\theta}^{(t)} - \bm{\theta}^{*}),\\
    &= \Vert\bm{\theta}^{(t)}\Vert_2^2 - 2 \bm{\theta}^{(t)T}\bm{\theta}^{*} + \Vert\bm{\theta}^{*}\Vert_2^2.
\end{split}
\label{Eq:C_t}
\end{equation}

\begin{equation}
\small
\begin{split}
    D^{(t+1)}
    &= (\bm{\theta}^{(t+1)} - \bm{\theta}^{*})^T(\bm{\theta}^{(t+1)} - \bm{\theta}^{*}),\\
    &= \Vert\bm{\theta}^{(t+1)}\Vert_2^2 - 2 \bm{\theta}^{(t+1)T}\bm{\theta}^{*} + \Vert\bm{\theta}^{*}\Vert_2^2, \\
    &= \Vert \bm{\theta}^{(t)} - \eta G^{(t)} \Vert_2^2 - 2(\bm{\theta}^{(t)} - \eta G^{(t)})^T\bm{\theta}^{*} + \Vert\bm{\theta}^{*}\Vert_2^2.
\end{split}
\label{Eq:C_t+1}
\end{equation}

According to Equation (\ref{Eq:C_t}) and (\ref{Eq:C_t+1}), the convergence rate can be simplified as follows.
\begin{equation}
\small
    \begin{split}
        C^{(t)} 
        &= - \mathbb{E}[D^{(t+1)} - D^{(t)}], \\
        &= - \mathbb{E}[\Vert \bm{\theta}^{(t)} - \eta G^{(t)} \Vert_2^2 + 2\eta G^{(t)T} \bm{\theta}^{*} - \Vert\bm{\theta}^{(t)}\Vert_2^2], \\
        &= - \mathbb{E}[-2\eta(\bm{\theta}^{(t)} - \bm{\theta}^{*})G^{(t)} + \eta^2G^{(t)T}G^{(t)}], \\
        &= 2\eta(\bm{\theta}^{(t)} - \bm{\theta}^{*})\mathbb{E}[G^{(t)}] - \eta^2 \mathbb{E}[G^{(t)}]^T\mathbb{E}[G^{(t)}], \\ 
        &- \eta^2\mathrm{Tr}(\mathbb{V}[G^{(t)}]), 
    \end{split}
\end{equation}
where $\mathrm{Tr}$ denotes trace, and $\mathbb{V}$ denotes variance. 

The expectation of the gradient is independent of the sampling probability $p_\mathrm{x}$.
maximizing convergence rate is equivalent to minimize $\mathrm{Tr}(\mathbb{V}[G^{(t)}])$, and \citet{needell2014stochastic, zhao2015stochastic} have demonstrated that the best sampling distribution is determined by $ q_{\mathrm{x}}^{(t)*} \propto \Vert \nabla_{\bm{\theta}^{(t)}} \ell(\mathrm{x}, \bm{\theta}^{(t)})\Vert_2 $ .

\subsection{Optimization Problem of PINNs}
As mentioned in the main content, PINNs learn approximate solutions $\hat{u}(t, \bm{x}; \bm{\theta})$ to fit the latent solutions $u(t, \bm{x})$ of PDEs.
\begin{equation}
    \begin{split}
        &\frac{\partial{u}}{\partial{t}} + \mathcal{N}_{\bm{x}} [u] = 0, \bm{x} \in \Omega, t\in[0, T], \\
        &u(t, \bm{x})|_{t=0} = u_0(\bm{x}), \bm{x} \in \Omega, \\
        &u(t, \bm{x})=g(t, \bm{x}), x\in \partial{\Omega}, t\in[0, T],
    \end{split}
    \label{Eq:AppendixProblemDefine}
\end{equation}
where $\bm{\theta}$ is the parameter of PINNs, $\mathcal{N}_{\bm{x}}$ denotes a differential operator consisted of spatial derivatives, $u_{0}(\bm{x})$ is the initial condition, $g(\bm{x})$ is the boundary condition, $\bm{x}$ is a D-dimensional position vector, and $\Omega$ is a subset of $\mathbb{R}^D$ with boundary $\partial{\Omega}$.

PINNs regard boundary conditions and initial conditions as penalty terms and formulate the constrained optimization problem into an unconstrained optimization problem:
\begin{equation}
    \bm{\theta^*} = \arg \min_{\bm{\theta}} r_f(\bm{\theta}) + \lambda_1 r_i(\bm{\theta}) + \lambda_2 r_b(\bm{\theta}),
    \label{Eq:AppendixOpProblemPINNs}
\end{equation}
where $\lambda_1$ and $\lambda_2$ are weights, $r_f(\bm{\theta})$, $r_i({\bm{\theta}})$ and $r_b({\bm{\theta}})$ are the residuals of PDEs, initial conditions, and boundary conditions, respectively.

The residual of PDEs is an integral over the spatial-temporal domain:
\begin{equation}
    r_f(\bm{\theta}) = \int_{\Omega \times [0, T]}(\frac{\partial{\hat{u}}}{\partial t} + \mathcal{N}_{\bm{x}}[\hat{u}])^2d\bm{x}dt. 
\end{equation}

The residual of initial conditions is an integral over the spatial domain with $t=0$:
\begin{equation}
    r_i(\bm{\theta}) = \int_{\Omega}(\hat{u}(0, \bm{x}) - u_0(\bm{x}))^2d\bm{x}.
\end{equation}

The residual of boundary conditions is an integral over the boundary of the domain:
\begin{equation}
    r_b(\bm{\theta}) = \int_{\partial\Omega \times [0, T]}(\hat{u}(t, \bm{x}) - g(t, \bm{x}))^2d\bm{x}dt.
\end{equation}

Equation \ref{Eq:AppendixOpProblemPINNs} can be approximated as a loss function of data points:
\begin{equation}
\begin{split}
    \bm{\theta}^* 
    &\approx  \arg\min_{\bm{\theta}}\mathcal{L}, \\ 
    & = \arg\min_{\bm{\theta}}\mathcal{L}_f + \lambda_1\mathcal{L}_b + \lambda_2 \mathcal{L}_i,
\end{split}
\label{Eq:AppendixOpProblemReal}
\end{equation}
where $\mathcal{L}$ is the total loss, $\mathcal{L}_f$, $\mathcal{L}_i$, and $\mathcal{L}_b$ are the losses of PDE residuals, initial conditions, and boundary conditions, respectively. The datasets of these losses are denoted by $N_f$, $N_i$, and $N_b$, respectively.

\begin{equation}
\begin{split}
    \mathcal{L}_f 
    &= \frac{1}{\vert N_f \vert}\sum_{\mathrm{x}\in N_f}\ell_f(\mathrm{x};\bm{\theta}), \\
    &= \frac{1}{\vert N_f \vert }\sum_{i=1}^{\vert N_f \vert} \Vert\frac{\partial{\hat{u}}}{\partial{t_i}} + \mathcal{N}_{\bm{x}_i}[\hat{u}(t_i, \bm{x_i})]\Vert_2^2,
\end{split}
\end{equation}
where $\ell_f$ is the loss of a data point $\mathrm{x} \in N_f$.

\begin{equation}
    \begin{split}
        \mathcal{L}_i 
        &= \frac{1}{\vert N_i \vert}\sum_{\rm{x}\in N_i}\ell_i(\rm{x};\bm{\theta}),\\
        &= \frac{1}{\vert N_i\vert}\sum_{i=1}^{\vert N_i\vert} \Vert\hat{u}(t_i, \bm{x}_i) - g(t_i, \bm{x}_i)\Vert^2_2,
    \end{split}
\end{equation}
where $\ell_i$ is the loss of a data point $\mathrm{x} \in N_i$.

\begin{equation}
\begin{split}
    \mathcal{L}_b 
    &= \frac{1}{\vert N_b \vert}\sum_{\rm{x}\in N_b}\ell_b(\rm{x};\bm{\theta}), \\ 
    &=\frac{1}{\vert N_b \vert}\sum_{i=1}^{\vert N_b \vert} \Vert\hat{u}(0, \bm{x}_i) - u_0(\bm{x}_i)\Vert^2_2,
\end{split}
\end{equation}
where $\ell_b$ is the loss of a data point $\mathrm{x} \in N_b$.

In DMIS, importance sampling is integrated into the training, and the detail of DMIS is described in the main content.

\section{Appendix C: Additional Experiment}
In this section, we supplement the experiments of solving the Diffusion Equation and the Allen-Cahn (AC) equation.

\subsection{Benchmark}
The same method has been used as described in the main content to divide the entire time domain $[0, T]$ into three segments: $[0, T/2]$, $[T/2, 3T/4]$, and $[3T/4, T]$. 
These three segments are used for training, validating, and testing.

\begin{table*}[!ht]
    \centering
    \begin{tabular*}{\hsize}{@{}@{\extracolsep{\fill}}c c c c c c c c c c@{}}
    \hline
    
    \hline
    \multicolumn{1}{c}{\multirow{2}{*}{PDE}} & \multicolumn{2}{c}{Network Config} & Optimizer Config & \multicolumn{3}{c}{DMIS Config} & \multicolumn{3}{c}{Dataset Config}\\ \cline{2-10} & Depth & Width & Learning rate & $\vert S \vert$ & $\gamma$ & $\beta$ & $\vert N_{f} \vert $ & $\vert N_{i} \vert $ & $\vert N_{b} \vert $ \\ \hline
        Diffusion & 4 & 32 & 0.002 & 1000 & 0.4 & 2 & 100000 & 2000 & 2000 \\  
        Allen-Cahn & 5 & 64 & 0.001 & 1000 & 0.4 & 1.5 & 60000 & 2000 & 2000\\
        \hline
        
        \hline
    \end{tabular*}
    \vspace{-0.2cm}
    \caption{The hyper-parameters used for each benchmark. $\vert S \vert$ is the set size of mesh points, $\gamma$ is the mesh update threshold, and $\beta$ is the hyper-parameter of reweighting. $\vert N_f \vert$, $\vert N_i \vert$, and $\vert N_b \vert$ are the dataset size of PDE residuals, initial conditions, and boundary conditions, respectively.}
    \label{Tb:Appendixhyperparameters}
\end{table*}

\begin{table*}[!ht]
    \centering
    \begin{tabular*}{\hsize}{@{}@{\extracolsep{\fill}}ccccccc@{}} 
    \hline
    
    \hline
    \multicolumn{1}{c}{\multirow{2}{*}{Method}} &
    \multicolumn{3}{c}{Diffusion Equation} &
    \multicolumn{3}{c}{Allen-Cahn Equation}  \\
    \cline{2-7} & ME & MAE & RMSE & ME & MAE & RMSE \\ \hline
    PINN-O & 0.181 & 0.089 & 0.095 & 1.406 & 0.133 & 0.232 \\ 
    PINN-N & 0.158 & 0.076 & 0.081 & 1.404 & 0.107 & 0.203 \\ \hline
    \textbf{PINN-DMIS(ours)} & \textbf{0.053} & \textbf{0.016} & \textbf{0.019} & \textbf{1.037} & \textbf{0.086} & \textbf{0.150} \\ 
    \hline
    
    \hline
    \end{tabular*}
    \vspace{-0.2cm}
    \caption{Comparison with PINN-O \citep{raissi2019physics} and PINN-N \citep{nabian2021efficient} on benchmarks of the Diffusion Equation and the Allen-Cahn Equation.}
    \label{Tb:AppendixAccEvaluation}
    \vspace{-10pt}
\end{table*}

\subsubsection{Diffusion Equation}
The diffusion equation is a crucial partial differential equation used to describe the change of matter density in the diffusion phenomenon. This equation can also describe the heat conduction process.
We consider the diffusion equation with an initial condition $u(0, x) = 2sin(\pi x) + 2x - 2x^3$ and a boundary condition $u(t, x) = 0, x \in \{0, 1\}$.
\begin{equation}
    \frac{\partial{u}}{\partial{t}}=1.2\frac{\partial^2{u}}{\partial{x^2}} + 5e^{-t}x, x\in [0, 1], t\in[0, 1].
\end{equation}

\subsubsection{Allen-Cahn Equation}
The Allen-Cahn (AC) equation is a type of reaction-diffusion equation. This equation is also the main mathematical model for phase-field simulation in material science, which describes the phase separation process of multi-component alloys. We consider the AC equation with an initial condition $u(0, x)=x^2 \cos (\pi x)$ and periodic boundary conditions:
\begin{equation}
\frac{\partial{u}}{\partial{t}} = \frac{0.001}{\pi}\frac{\partial^2{u}}{\partial{x^2}} + 2\sin(\pi u), x\in[-1, 1], t\in[0, 1].
\end{equation}

\subsection{Experimental Setting}

\subsubsection{Hyper-parameter}
Table \ref{Tb:Appendixhyperparameters} summarizes hyper-parameters for each additional benchmark.
The Adam optimizer is used on all additional benchmarks. 
For PINN-N, we uniformly sample 10,000 seeds within the domain, and other parameters are consistent with Table \ref{Tb:Appendixhyperparameters}.
Experiments on each additional benchmark have the same random seed.

\subsubsection{Evaluation Metrics}
We use the same evaluation metrics as described in the main content. 
The maximum error(AE), mean absolute error (MAE), and root mean square error (RMSE) are used to evaluate prediction accuracy.
The calculation time $TC$ and the number of iteration steps $NC$ required for convergence are used to evaluate convergence behavior.

\subsection{Experimental Result}

\subsubsection{Diffusion Equation}
\begin{figure}[!ht]
    \centering
    \includegraphics[width=0.88\columnwidth]{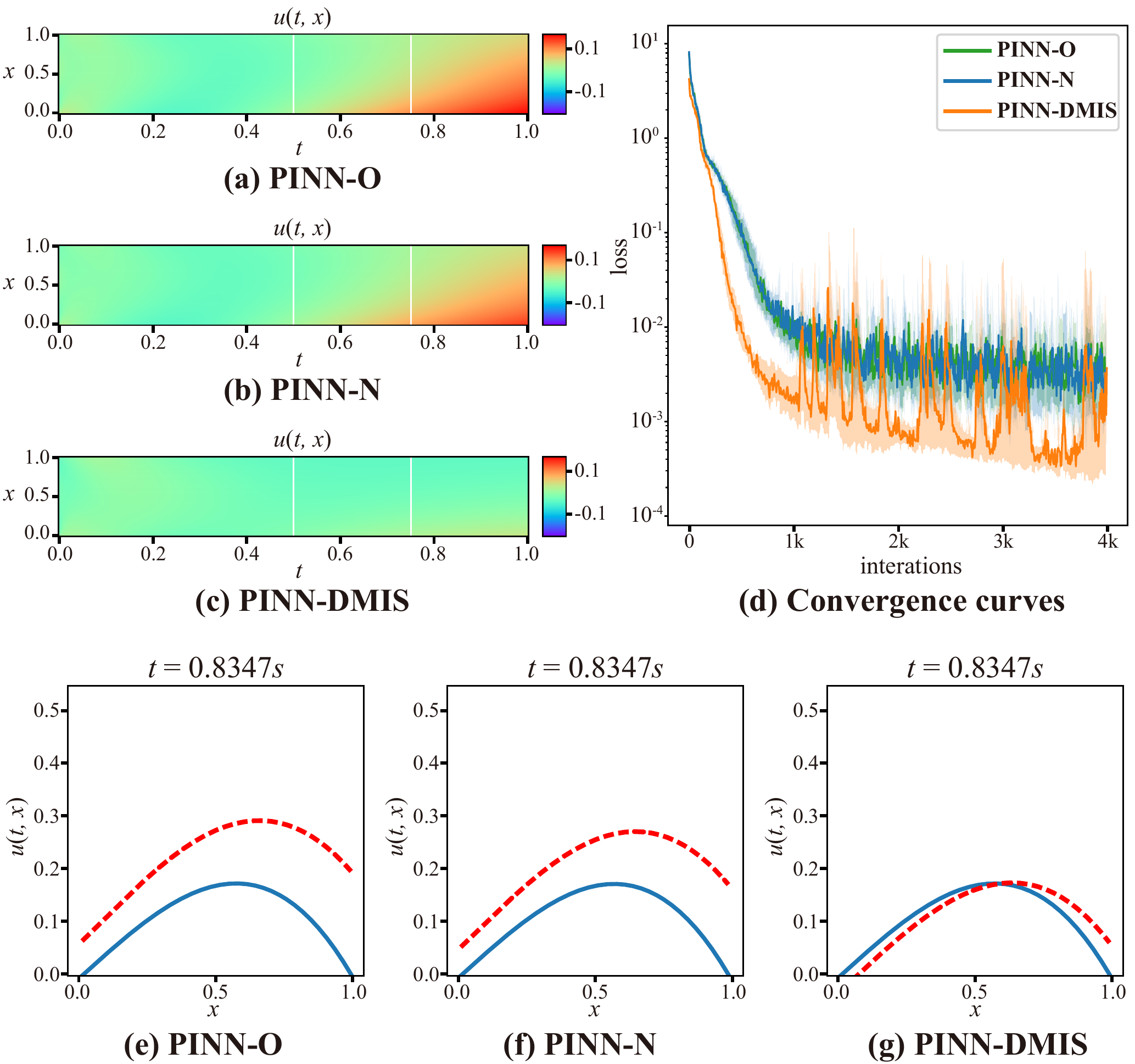}
    \vspace{-0.2cm}
    \caption{The Diffusion Equation. (a, b, c) The prediction errors of PINN-O, PINN-N, and PINN-DMIS, respectively. (d) Convergence curves. (e, f, g) show the prediction (red) and ground truth (blue).}
    \label{fig:Diffusion}
\end{figure}

Figure \ref{fig:Diffusion}(a, b, c) report the prediction error by PINN-O, PINN-N, and PINN-DMIS for solving the diffusion equation. 
Compared with other state-of-the-art methods, DMIS has the best performance.
The snapshots in Figure \ref{fig:Diffusion} show that DMIS has the best prediction results in the whole domain.
Table \ref{Tb:AppendixAccEvaluation} summarizes the evaluation results of prediction accuracy for solving Diffusion Equation.
DMIS has the best performance on all evaluation metrics.
Figure \ref{fig:Diffusion}(d) reports the convergence curves for solving the Diffusion Equation.
DMIS converges fastest and has the lowest training error.
However, the loss value fluctuates greatly in DMIS. 
We speculate that the loss value fluctuates greatly due to the unstable sample weights. Affected by the complex initial conditions, the triangular mesh updates frequently during training.
Table \ref{Tb:DiffusionConvergenceBehavior} compares the convergence behaviors.
DMIS has the fastest convergence speed.

\subsubsection{Allen-Cahn Equation}
\begin{figure}[!ht]
    \centering
    \includegraphics[width=0.88\columnwidth]{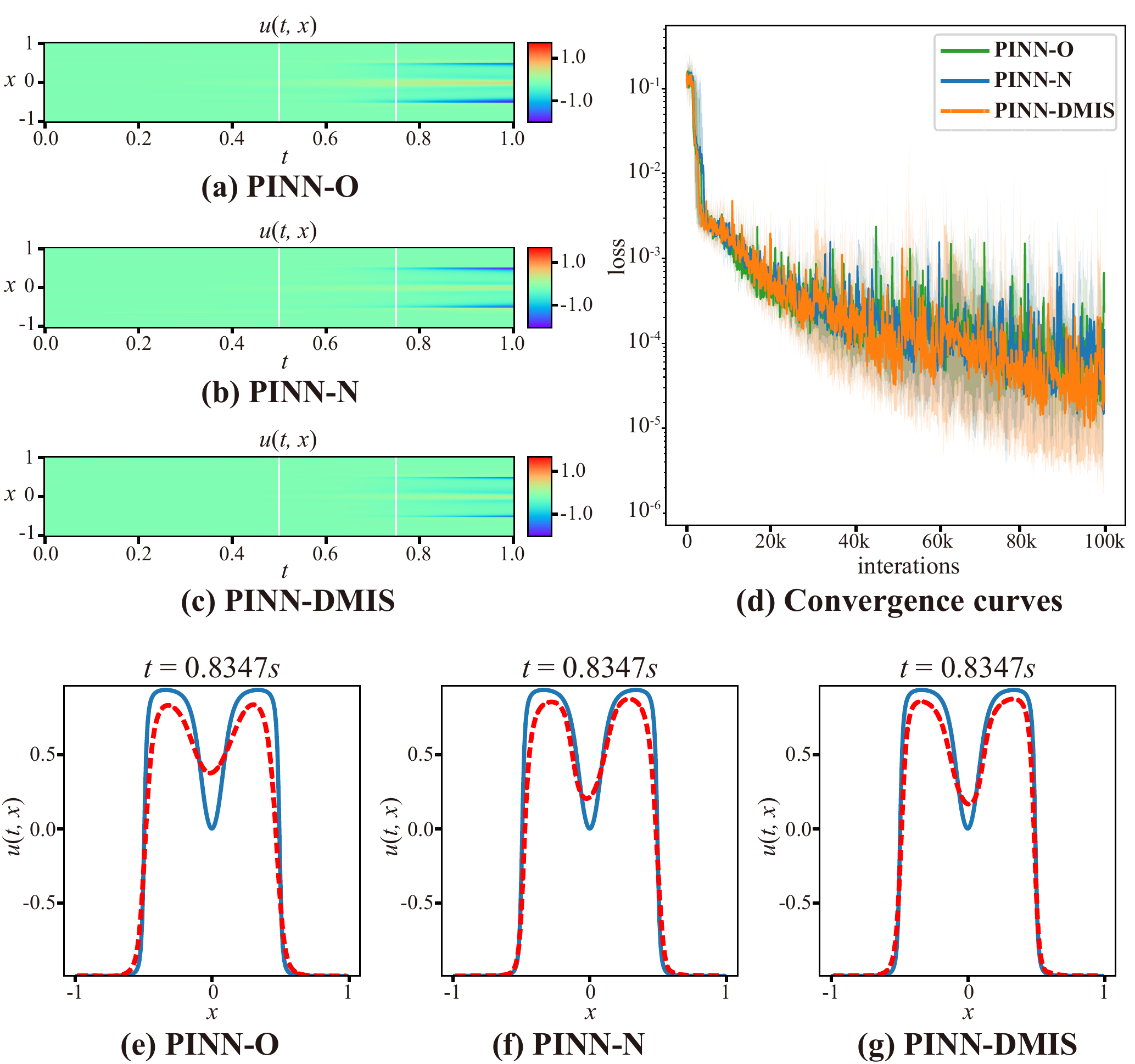}
    \vspace{-0.2cm}
    \caption{The Allen-Cahn Equation. (a, b, c) The prediction errors of PINN-O, PINN-N, and PINN-DMIS, respectively. (d) Convergence curves. (e, f, g) show the prediction (red) and ground truth (blue).}
    \label{fig:ACEq}
\end{figure}

Figure \ref{fig:ACEq}(a, b, c) report the prediction errors for solving the Allen-Cahn Equation. 
DMIS has the lowest prediction error.
Figure \ref{fig:ACEq}(e, f, g) report the snapshots of prediction at t=0.8347s.
Compared with PINN-O and PINN-N, DMIS has better prediction results on the position with a large gradient.
Table \ref{Tb:AppendixAccEvaluation} summarizes the evaluation results for solving the Allen-Cahn Equation.
DMIS has the best performance on all evaluation metrics.
Figure \ref{fig:ACEq}(d) reports the convergence curves for solving the Allen-Cahn Equation.
DMIS does not show obvious advantages.
The result shows that DMWE fails to generate an effective mesh.
The failure reasons are analyzed in detail in the next section.
Table \ref{Tb:ACEqConvergenceBehavior} reports the convergence behaviors for solving the Allen-Cahn Equation.
These three methods have similar convergence behavior and DMIS has a slender lead convergence speed. 

\begin{table}[!t]
    \centering
    \begin{tabular*}{\hsize}{@{}@{\extracolsep{\fill}}c c c c c@{}}
    \hline
    
    \hline
    ~ & $TC_1/s$ & $TC_2/s$ & $NC_1$ & $NC_2$ \\ \hline
    PINN-O & 19 & 42 & 1089 & 2420 \\
    PINN-N & 82 & 173 & 1089 & 2299 \\ \hline
    \textbf{PINN-DMIS(ours)} & \textbf{18} & \textbf{21} & \textbf{847} & \textbf{968} \\
    \hline
    
    \hline
    \end{tabular*}
    \vspace{-0.2cm}
    \caption{Evaluation results of convergence behavior for solving the Diffusion Equation}
    \label{Tb:DiffusionConvergenceBehavior}
\end{table}

\section{Appendix D: Limitations of DMIS}
Because DMIS is based on importance sampling, when importance sampling fails, DMIS also fails to speed up the training of PINNs. If all data points have the same loss value, the sampling probabilities of data points are equal to $1/N_f$, and the importance sampling degenerates into uniform sampling. In this case, DMIS is not only unable to accelerate training, but also causes additional computational costs.

In practice, considering the loss values of data points are not precisely equal, DMIS is always valid. However, we still need to pay attention to the balance between the improvement and the additional calculation cost brought by DMIS.
In our experiments, DMIS shows excellent performance in solving the Schrödinger equation. In solving the Schrödinger equation, DMIS can significantly improve the convergence efficiency because the loss values of data points are quite different. On the other hand, in solving the Allen-Cahn equation, DMIS almost degenerates into uniform sampling.

\begin{table}[!t]
    \centering
    \begin{tabular*}{\hsize}{@{}@{\extracolsep{\fill}}c c c c c@{}}
    \hline
    
    \hline
    ~ & $TC_1/s$ & $TC_2/s$ & $NC_1$ & $NC_2$ \\ \hline
    PINN-O & 70 & \textbf{97} & 2847 & 3989 \\
    PINN-N & 213 & 304 & 3181 & 4562 \\ \hline
    \textbf{PINN-DMIS(ours)} & \textbf{67} & 104 & \textbf{2243} & \textbf{3459} \\
    \hline
    
    \hline
    \end{tabular*}
    \vspace{-0.2cm}
    \caption{Evaluation results of convergence behavior for solving the Allen-Cahn Equation}
    \label{Tb:ACEqConvergenceBehavior}
\end{table}

In addition to the limitation brought by importance sampling, the triangular mesh also affects the effect of DMIS. 
For DMIS, the instability of the loss value distribution leads to frequent updating of the triangular mesh, which affects the training efficiency and stability of training. 
In solving the diffusion equation, the mesh updates frequently, and the loss value fluctuates greatly, as shown in Figure \ref{fig:Diffusion}(d).
These defects can be alleviated by setting an appropriate updating threshold. 
The selection of the threshold mainly depends on experience. Fortunately, our experiments show that the impact of the threshold is relatively stable within a reasonable range $[0.3, 0.6]$ and we recommend setting the threshold $\gamma$ to 0.4.

In this paper, we only verify the excellent performance and generalization ability of DMIS on five widely-used benchmarks. We will verify the applicability of DMIS on more complex PDEs, for example, high-dimensional PDEs, in the future.

\end{document}